\documentclass[11pt]{article}
\usepackage{amssymb}
\usepackage{color, amsmath,amssymb, amsfonts, amstext,amsthm, latexsym}

\usepackage{amssymb, epsfig, amssymb, latexsym}
\usepackage{amsmath}
\usepackage{graphicx}
\usepackage{longtable}

\newtheorem{theorem}{Theorem}
\newtheorem{lemma}{Lemma}

  \newtheorem{Main Results}[theorem]{MainResults}

 \newtheorem{remark}[lemma]{Remark}

\newcommand{\PX}{{\Bbb{P}}}

\newcommand{\Om}{\Omega}

\begin{document}

\title{Fokker-Planck equations for nonlinear dynamical systems driven by non-Gaussian L\'evy processes
\footnote{This work was supported by NSFC grants 10971125 and 11028102, the NSF grant DMS-1025422, and the fundamental research funds for the Central Universities HUST 2010ZD037. } }

 \author{Xu Sun\\
 School of Mathematics and Statistics\\Huazhong University of Science and Technology\\ Wuhan 430074, Hubei, China \\
   E-mail: xsun15@gmail.com\\ \\
Jinqiao Duan\\
Institute for Pure and Applied Mathematics, University of California\\
Los Angeles, CA 90095\\ \& \\
Department of Applied Mathematics, Illinois Institute of Technology \\
   Chicago, IL 60616, USA \\
   E-mail: duan@iit.edu\\
   \\
 }

\date{January 9, 2012}
\maketitle

\newpage
\begin{abstract}
The Fokker-Planck equations describe time evolution of probability densities of stochastic dynamical systems and are thus widely used to quantify random phenomena such as uncertainty propagation.  For dynamical systems driven by non-Gaussian L\'evy processes, however, it is   difficult to obtain explicit forms of Fokker-Planck equations because the adjoint operators of the associated infinitesimal generators usually do not have exact formulation. In the present paper,    Fokker-Planck equations are derived in terms of infinite series for nonlinear stochastic differential equations with non-Gaussian L\'evy processes. A few examples are presented to illustrate the method.

\bigskip

PACS 2010:  05.40.Ca;  05.10.Gg;  02.50.Fz


\end{abstract}

\maketitle

\pagestyle{plain}

\section{Introduction and statement of the problem}  \label{intro}

The Fokker-Planck equations are one of the  deterministic   tools to  quantify how   randomness propagates or  evolves in nonlinear dynamical systems. For   stochastic differential equations (SDEs) with Gaussian processes such as Brownian motion,      the Fokker-Planck equations  are well established  \cite{Okse2003, Klebaner2005}. However, for SDEs with  non-Gaussian processes such as L\'{e}vy processes, explicit forms of the Fokker-Planck equations are not easily available, except in some special cases \cite{Applebaum2009, Schertzer}.
The difficulty is to obtain the expressions for the adjoint operators of the infinitesimal generators  associated with these SDEs.

  L\'{e}vy processes are a   class of  stochastic processes having independent  and stationary increments, as well as  stochastically continuous sample paths \cite{Applebaum2009, Sato1999}. Given a sample space $ \Om$, together with a probability measure $ \PX$ and the corresponding mathematical expectation $\mathbf{E}$.  A  L\'{e}vy process $L_t$, taking values in $\mathbb{R}^d$, is
characterized by a drift parameter $b\in\mathbb{R}^d$ , a
$d\times d$ positive-definite covariance  matrix $A$ and a    measure $\nu$ defined on
$\mathbb{R}^d$ and concentrated on $\mathbb{R}^d\backslash\{0\}$. In fact, this measure $\nu$
  satisfies the following condition \cite{Applebaum2009}
\begin{align}\label{s1_1}
\int_{\mathbb{R}^d\backslash \{0\} }(\|y\|^2\wedge 1)\nu(dy)<\infty,
\end{align}
or equivalently
\begin{align}\label{s1_2}
\int_{\mathbb{R}^d\backslash \{0\}}\frac{\|y\|^2}{1+\|y\|^2}\nu(dy)<\infty.
\end{align}
Here $\|\cdot \|$ is the usual Euclidean norm or length in $\mathbb{R}^d$.  This measure $\nu$ is  called a L\'{e}vy jump measure for  the
L\'{e}vy process $L_t$. A L\'{e}vy process with the generating
triplet $(b, A, \nu)$ has the L\'{e}vy-It\^{o} decomposition
\begin{align}\label{s1_3}
L_t=bt+B_t+\int_{\|y\|< 1}y\tilde{N}(t, dy)+\int_{\|y\|\geq 1}y{N}(t,
dy),
\end{align}
where $N(dt, dx)$ is the Poisson random measure, $\tilde{N}(dt,
dx)=N(dt, dx)-\nu(dx)dt$  is the compensated Poisson random
measure, and $B_t$ is an independent $d$-dimensional Brownian motion (i.e., Wiener process) with
covariance matrix $A$.
Equation (\ref{s1_3}) can be formally rewritten in a differential form as
\begin{align}\label{s1_4}
{\rm d} L_t =b{\rm d} t+{\rm d} B_t+ \int_{\|y\|< 1}y\tilde{N}({\rm d} t, {\rm d} y)+\int_{\|y\|\geq 1}y{N}({\rm d} t,
{\rm d} y).
\end{align}

 We shall consider stochastic dynamical systems described by the following SDE   in the It\^{o} form
\begin{equation}\label{sde_Ito}
dX_t=f(X_t, t)dt+\sigma(X_{t-}, t)  dL_t, \;\;   X_0=x,
\end{equation}
or   in the Marcus form
\begin{equation}\label{Mainequation_strat}
{\rm d}X_t=f(X_t, t){\rm d}t+\sigma(X_{t-}, t)\diamond {\rm d}L_t, \;\;  X_0=x,
\end{equation}
where $L_t$ is a  L\'{e}vy process with the generating triplet     $(b, A, \nu)$.  Equation (\ref{sde_Ito}) is meaningful in the sense of
\begin{equation}\label{sde_Ito_1}
X_t=X_0 + \int_0^t f(X_s, s)ds+ \int_0^t \sigma(X_{s-}, s)  dL_s,
\end{equation}
where the last term is an It\^o integral. Note that L\'{e}vy processes are semimartingales, and It\^o integrals with respect to semimartingales are well defined \cite{Protter2004}. The  equation  (\ref{Mainequation_strat})    is interpreted as
\begin{align}\label{Mainequation_strat_1}
X_t =X_0 + \int_0^t f(X_s, s){\rm d}s+ \int_0^t \sigma(X_{s-}, s)\diamond {\rm d} L_s,
 \end{align}
where      ``$\diamond$" indicates Marcus integral \cite{Marcus1978, Marcus1981, Applebaum2009}   defined by
\small
\begin{align}\label{section1_tmp1}
 \int_0^t \sigma(X_{s-}, s)\diamond {\rm d} L_s&=
\int_0^t \sigma(X_{s-}, s)  {\rm d}L_s + \frac{1}{2} \int_0^t \sigma(X_{s-},s) \sigma^\prime (X_{s-},s) {\rm d} [L_s, L_s]^c\nonumber\\
 &\quad + \sum_{0\leq s \leq t} \left[ \xi(\Delta L(s), \sigma(X_{s-}, s),  X_{s-}) -X_{s-} -\sigma(X_{s-},s)\Delta L_{s} \right],
 \end{align}
 \normalsize
 with $ \xi(r, g(x),  x)$   being the value at $z=1$ of the solution of the following ordinary differential equation:
\begin{align}
\frac{\rm d}{{\rm d} z} y(z) =rg(y(z)), \quad\quad y(0)=x.
\end{align}

\textbf{Assumption ($H_1$):}\\
We assume that the appropriate Lipschitz and growth conditions on the  drift $f$ and noise intensity $\sigma$ are satisfied so that both  equations  \eqref{sde_Ito} and \eqref{Mainequation_strat}  have unique solutions.

\bigskip

 The main objective of this paper is to derive an expression of the Fokker-Planck equations for nonlinear dynamical systems described by SDEs  (\ref{sde_Ito}) and (\ref{Mainequation_strat}), in \S \ref{ito888} and \S \ref{marcus888}, respectively.   For simplicity, we only consider one-dimensional case, i.e., $X_t$ and $L_t$ in (\ref{sde_Ito}) and (\ref{Mainequation_strat}) are all scalar processes. The conclusion can be generalized into higher dimensional cases which, however, will not be considered in this paper.
A few examples will be presented in \S \ref{examples888}.

\section{Fokker-Planck equations for SDEs with It\^o integrals}
\label{ito888}

 In this section, we derive the Fokker-Planck equation for It\^o SDE (\ref{sde_Ito}),  in which $L_t$ being a scalar L\'{e}vy process  with the triplet  $(b, A, \nu)$. In this case, $A$ is a non-negative scalar.

 Substituting (\ref{s1_4}) into (\ref{sde_Ito}), we get
\begin{align}\label{s31_1}
dX_t &=f(X_t, t) {\rm d} t+b\sigma(X_{t-}, t) {\rm d}t + \sigma (X_{t-},t) {\rm d} B_t + \int_{|y|< 1} \sigma (X_{t-}, t) y\,\tilde N({\rm d}t,{\rm d}y )\nonumber\\
&\quad + \int_{|y| \ge 1} \sigma (X_{t-}, t) y\,N({\rm d}t,{\rm d}y ).
\end{align}
By It\^o formula \cite{Applebaum2009}, it follows from (\ref{s31_1}) that for any smooth function $\phi (x)$,
\footnotesize
\begin{align}\label{s31_2}
\phi (X_{t+\Delta t})-\phi (X_t) &= \int_t^{t+\Delta t} f(X_{s-}, s) \frac{\partial}{\partial x} \phi (X_s)\,{\rm d}s  +  \int_t^{t+\Delta t} b\sigma (X_{s-}, s) \frac{\partial}{\partial x} \phi (X_s)\,{\rm d}s\nonumber\\
 &\quad  + \int_t^{t+\Delta t}  \sigma (X_{s-}, s) \frac{\partial}{\partial x} \phi (X_s) \, {\rm d} B_s + \frac{A}{2} \int_t^{t+\Delta t} \sigma ^2(X_{s-},s)\frac{\partial ^2 }{\partial x^2 } \phi (X_s)\,{\rm d}s \nonumber\\
&\quad +\int_t^{t+\Delta t}\int_{|y|\ge 1} \left[ \phi\left(X_{s-} + y\sigma (X_{s-},s)\right)-\phi(X_{s-})\right]\,N({\rm d}s, {\rm d} y)\nonumber\\
&\quad + \int_t^{t+\Delta t}\int_{|y|< 1} \left[ \phi\left(X_{s-} + y\sigma (X_{s-},s)\right)-\phi(X_{s-})\right]\,\widetilde N({\rm d}s, {\rm d} y)\nonumber\\
&\quad + \int_t^{t+\Delta t}\int_{|y|< 1} \left[ \phi\left(X_{s-} + y\sigma (X_{s-},s)\right)-\phi(X_{s-})- y\sigma (X_{s-}, s) \frac{\partial }{\partial x} \phi(X_{s-}) \right]\,\nu ( {\rm d} y){\rm d}s\,.
\end{align}
\normalsize
Let $Q$ be the infinitesimal generator associated with the Markovian solution process $X_t$. Then for any function $\phi$ in $C_0^\infty (\mathbb R)$, consisting of smooth functions with compact support  on $\mathbb{R}$,  we obtain
\begin{align}\label{s31_3}
Q\phi(x)&=\lim_{\Delta t\to 0} \frac{\mathbf E\{ \phi (X_{t+\Delta t})\big|X_t=x)\}-\phi (x)}{\Delta t}\nonumber\\
&=f(x,t)\frac{\partial }{\partial x} \phi (x) +b\, \sigma(x,t)\frac{\partial }{\partial x} \phi (x) + \frac{A}{2} \sigma^2(x,t) \frac{\partial^2}{\partial x^2} \phi(x)\nonumber\\
&\quad +  \int_{|y|< 1} \left[ \phi(x+y\sigma (x,t)) -\phi(x)-  y \sigma (x,t) \frac{\partial}{\partial x} \phi(x)\right]\,\nu ({\rm d} y)\nonumber\\
&\quad +  \int_{|y| \ge 1} \left[ \phi(x+y\sigma (x,t)) -\phi(x) \right]\,\nu ({\rm d} y),
\end{align}
or equivalently,
\begin{align}\label{tt_3}
Q\phi(x)&=f(x,t)\frac{\partial }{\partial x} \phi (x) +b\, \sigma(x,t)\frac{\partial }{\partial x} \phi (x) + \frac{A}{2} \sigma^2(x,t) \frac{\partial^2}{\partial x^2} \phi(x)\nonumber\\
&\quad +  \int\limits_{R\backslash \{0\} }\left[ \phi(x+y\sigma (x,t)) -\phi(x)-I_{(-1, 1)}(y) y \sigma (x,t) \frac{\partial}{\partial x} \phi(x)\right]\,\nu ({\rm d} y),
\end{align}
where $I_{(-1, 1)}$ the indicator function of the set $(-1, 1)$.
To get the second identity in (\ref{s31_3}), we have used the fact that
\begin{align}\label{s31_3a}
\mathbf E\left\{\int_t^{t+\Delta t}  \sigma (X_{s-}, s) \frac{\partial}{\partial x} \phi (X_s) \, {\rm d} B_t\Bigg | X_t=x
\right\}=0,
\end{align}
\begin{align}\label{s31_3b}
\mathbf E\left\{\int_t^{t+\Delta t}\int_{|y|< 1} \left[ \phi\left(X_{s-} + y\sigma (X_{s-})\right)-\phi(X_{s-},s)\right]\,\widetilde N({\rm d}s, {\rm d} y) \Bigg | X_t=x\right\}=0,
\end{align}
and
\begin{align}
&\mathbf E\left\{\int_t^{t+\Delta t}\int_{|y|\ge 1} \left[ \phi\left(X_{s-} + y\sigma (X_{s-})\right)-\phi(X_{s-},s)\right]\,N({\rm d}s, {\rm d} y) \Bigg | X_t=x\right\}\nonumber\\
&\quad =  \int_t^{t+\Delta t}\int_{|y|\ge 1} \left[ \phi\left(x + y\sigma (x,s)\right)-\phi(x)\right]\,  \nu ({\rm d} y){\rm d}s  .
\end{align}
Denoting $u(x,t)=\mathbf{E}\left(\phi(X_t)\mid X_0=x\right)$, it can then be shown that \cite{Applebaum2009}
\begin{align}\label{s31_11}
\frac{\partial}{\partial t}u(x,t) =Qu(x,t),
\end{align}
which is the well known backward Kolmogorove equation. Let $p(y,t)$ be the  probability density function for $X_t$ associated with the SDE (\ref{sde_Ito}). Then  (\ref{s31_11}) becomes
\begin{align}\label{tt_1}
\int_\mathbb{R} \frac{\partial}{\partial t}\left( \phi(y) p(y, t)\right)\, {\rm d}y = \int_{\mathbb{R}} Q \phi(y) \, p(y,t) \,{\rm d}y.
\end{align}
Let $Q^*$ be the adjoint operator of $Q$.  Therefore,
\begin{align}\label{december_13_tt1}
\int_{\mathbb{R}} Q\phi(y) \, p(y,t) \,{\rm d}y = \int_{\mathbb{R}}  \phi(y) Q^* p(y,t) \,{\rm d}y.
\end{align}
With relation (\ref{december_13_tt1}) in mind, equation (\ref{tt_1}) becomes
\begin{align}\label{tt_2}
\int_\mathbb{R} \phi(y)  \left(\frac{\partial}{\partial t} p(y,t)- Q^* p(y,t)\right)\, {\rm d}y =0.
\end{align}
Since (\ref{tt_2}) is true for any $\phi\in C_0^\infty (\mathbb{R})$, we get
\begin{align}\label{tt_5}
\frac{\partial}{\partial t} p(y,t)= Q^* p(y,t),
\end{align}
which is the Fokker-Planck equation,  the governing equation for the transition probability density $p$ for the stochastic dynamical system  (\ref{sde_Ito}).

\bigskip

Now we try to derive an expression for the adjoint operator $Q^*$.
It follows from (\ref{tt_3}) that the operator $Q$ can be written as
\begin{align}
Q=A_1+ A_2
\end{align}
where $A_1$ and $A_2$ are defined as
\begin{align}
A_1 \phi (x) = f(x,t)\frac{\partial }{\partial x} \phi (x) +b\, \sigma(x,t)\frac{\partial }{\partial x} \phi (x) + \frac{A}{2} \sigma^2(x,t) \frac{\partial^2}{\partial x^2} \phi(x),
\end{align}
and
\begin{align}\label{ttt_1}
A_2 \phi(x) =\int\limits_{R\backslash \{0\} }\left[ \phi(x+y\sigma (x,t)) -\phi(x)-I_{(-1, 1)}(y) y \sigma (x,t) \frac{\partial}{\partial x} \phi(x)\right]\,\nu ({\rm d} y),
\end{align}
 respectively.
Note that $A_1^*$ is expressed as
 \begin{align}\label{dec13_tt3}
 A_1^* p(x,t)= - \frac{\partial}{\partial x}\left[ p(x,t)\left(f(x,t) + b \sigma(x,t)\right)\right] + \frac{A}{2} \frac{\partial^2}{\partial x^2} ( \sigma^2(x,t)p(x,t)).
 \end{align}
We now find an expression for $A_2^*$. Once $A_2^*$ is obtained, the adjoint operator $Q^*$ can be expressed as
\begin{align}\label{dec13_tt4}
Q^*= A_1^* + A_2^*
\end{align}
and then the Fokker-Planck equation $p_t=Q^*p$ will be obtained.

\medskip

 By Taylor expansion, $ \phi(x+y\sigma (x))$ can be written as
\begin{align}\label{s31_4}
 \phi(x+y\sigma (x,t)) =  \phi(x ) + \sum_{k=1}^{\infty} \frac{y^k}{k!}\sigma^k (x, t) \frac{\partial ^k}{\partial x^k}\phi (x).
 \end{align}
 Substitute (\ref{s31_4}) into (\ref{tt_3}) and (\ref{ttt_1}), respectively,  we get
 \begin{align}\label{s31_5a}
Q\phi(x)&=f(x,t)\frac{\partial }{\partial x} \phi (x) +b\, \sigma(x,t)\frac{\partial }{\partial x} \phi (x) + \frac{A}{2} \sigma^2(x,t) \frac{\partial^2}{\partial x^2} \phi(x)\nonumber\\
 &\quad + \int\limits_{R\backslash \{0\} }\left[ \sum_{k=1}^{\infty} \frac{y^k}{k!}\sigma^k (x,t) \frac{\partial ^k}{\partial x^k}\phi (x)-I_{(-1, 1)}(y) y \sigma (x,t) \frac{\partial}{\partial x} \phi(x)\right]\,\nu ({\rm d} y),
\end{align}
and
 \begin{align}\label{s31_5}
A_2\,\phi(x)=\int\limits_{R\backslash \{0\} }\left[ \sum_{k=1}^{\infty} \frac{y^k}{k!}\sigma^k (x,t) \frac{\partial ^k}{\partial x^k}\phi (x)-I_{(-1, 1)}(y) y \sigma (x,t) \frac{\partial}{\partial x} \phi(x)\right]\,\nu ({\rm d} y).
\end{align}

\medskip

It follows from (\ref{s31_5}) that
\small
\begin{align}\label{tt_7}
&\int_{\infty}^{\infty} A_2 \phi(x) p(x,t)\,{\rm d}x \nonumber\\
&\quad =\int\limits_{\mathbb{R}}\left(\int\limits_{\mathbb{R}\backslash \{0\} }\left[ \sum_{k=1}^{\infty} \frac{y^k}{k!}\sigma^k (x,t) \frac{\partial ^k}{\partial x^k}\phi (x)-I_{(-1, 1)}(y) y \sigma (x,t) \frac{\partial}{\partial x} \phi(x)\right] \, \nu( {\rm d} y) \right)p(x,t)\,{\rm d}x\nonumber\\
&\quad =\int\limits_{\mathbb{R}\backslash \{0\} }  \left( \int\limits_{\mathbb{R}}  \left[ \sum_{k=1}^{\infty} \frac{y^k}{k!}\sigma^k (x,t) \frac{\partial ^k}{\partial x^k}\phi (x)-I_{(-1, 1)}(y) y \sigma (x,t) \frac{\partial}{\partial x} \phi(x)\right]  p(x,t)\,{\rm d}x \right)\, \nu( {\rm d} y)\nonumber\\
&\quad =\int\limits_{\mathbb{R}\backslash \{0\} }  \left( \int\limits_{\mathbb{R}}  \left[ \sum_{k=1}^{\infty} \frac{(-y)^k}{k!} \frac{\partial ^k}{\partial x^k} \left(\sigma^k (x,t) p(x,t)\right) + I_{(-1, 1)}(y) y  \frac{\partial}{\partial x} \left(\sigma (x,t)p(x,t)\right) \right]  \phi(x)\,{\rm d}x \right)\, \nu( {\rm d} y)\nonumber\\
&\quad =  \int\limits_{\mathbb{R}} \left(\int\limits_{\mathbb{R}\backslash \{0\} }  \left[ \sum_{k=1}^{\infty} \frac{(-y)^k}{k!} \frac{\partial ^k}{\partial x^k} \left(\sigma^k (x,t) p(x,t)\right) + I_{(-1, 1)}(y) y  \frac{\partial}{\partial x} \left(\sigma (x,t)p(x,t)\right) \right]  \, \nu( {\rm d} y)\right) \phi(x)\,{\rm d}x\nonumber\\
&\quad = \int\limits_{\mathbb{R}}\phi (x, t) A_2^*\,p(x,t)  \,{\rm d}x,
\end{align}
\normalsize
where $A_2^*$, the adjoint of $A_2$, is  as follows
\begin{align}\label{dec13_tt2}
A_2^*\,p(x,t) =\int\limits_{\mathbb{R}\backslash \{0\} }  \left[ \sum_{k=1}^{\infty} \frac{(-y)^k}{k!} \frac{\partial ^k}{\partial x^k} \left(\sigma^k (x,t) p(x,t)\right) + I_{(-1, 1)}(y) y  \frac{\partial}{\partial x} \left(\sigma (x,t)p(x,t)\right) \right]  \, \nu( {\rm d} y).
\end{align}
Note that  in obtaining the third identity in (\ref{tt_7}), we have made use of
\begin{align}\label{dec13_tt5}
  \int\limits_{\mathbb{R}} p(x,t) \sigma^k (x,t) \frac{\partial ^k}{\partial x^k}\phi (x)\, {\rm d}x  = (-1)^k \int\limits_{\mathbb{R}} \phi(x) \frac{\partial ^k}{\partial x^k} \left(\sigma^k (x,t) p(x,t)\right) \,{\rm d}x, \quad\forall k\in\mathbb{N}.
\end{align}
Here $\mathbb{N}$ is the set of natural numbers.
Given  (\ref{dec13_tt2}), (\ref{dec13_tt3})  and (\ref{dec13_tt4}), it follows from (\ref{tt_5}) to get the desired Fokker-Planck equation. We summarize this result in the following theorem.

\begin{theorem} \label{fpe-ito}
Under the \textbf{Assumption ($H_1$)}, the   Fokker-Planck equation for the It\^o SDE (\ref{sde_Ito}) is
\small
\begin{align}\label{multiplicative_FPK}
\frac{\partial p(x,t)}{\partial t}
 &=-\frac{\partial }{\partial x}\left( f(x,t)p (x, t) \right) -b\frac{\partial }{\partial x}(\sigma(x,t)  p (x,t) ) + \frac{1}{2} A \frac{\partial^2}{\partial x^2} (\sigma^2(x,t) p(x,t)) \nonumber\\
 &\quad +\int\limits_{\mathbb{R}\backslash \{0\} }  \left[ \sum_{k=1}^{\infty} \frac{(-y)^k}{k!} \frac{\partial ^k}{\partial x^k} \left(\sigma^k (x,t) p(x,t)\right) + I_{(-1, 1)}(y) y  \frac{\partial}{\partial x} \left(\sigma (x,t)p(x,t)\right) \right]  \, \nu( {\rm d} y).
\end{align}
\normalsize
\end{theorem}

\begin{remark}
Two special cases are noted here.

(i)
When $\sigma(x,t)=1$, the Fokker-Planck equation  (\ref{multiplicative_FPK}) reduces to
\begin{align}\label{additive_FPK}
\frac{\partial p(x,t)}{\partial t}
 &=-\frac{\partial }{\partial x}\left( f(x,t)p (x, t) \right) -b\frac{\partial }{\partial x}  p (x,t) ) + \frac{A}{2}  \frac{\partial^2}{\partial x^2} \ p(x,t) \nonumber\\
 &\quad + \int\limits_{R\backslash \{0\} }\left( p(x-y,t ) -p(x,t)+I_{(-1, 1)}(y)y \frac{\partial}{\partial x}  p(x,t)  \right)\nu ({\rm d} y).
\end{align}
This is the Fokker-Planck equation for an It\^o SDE   with an additive L\'{e}vy process.

(ii)
When the L\'{e}vy jump measure $\nu$ satisfies $\int _{|x|\ge 1} |x|^k \,\nu ({\rm d} x) \le \infty$ for $k\in \mathbb{N}$, which is true if and only if $\mathbf{E} (|X_t|^k) \le \infty$ (see \cite{Applebaum2009}), the Fokker-Planck equation (\ref{multiplicative_FPK}) can be rewritten as
\begin{align}\label{multiplicative_FPK_1}
\frac{\partial p(x,t)}{\partial t}=-\frac{\partial }{\partial x}\left( f(x,t)p (x, t) \right) + \sum_{k=1}^{\infty}   c_k  \frac{\partial ^k}{\partial x^k} \left(\sigma^k (x,t) p(x,t)\right),
\end{align}
where
\begin{align}
c_1= -b - \int_{|y|\ge 1} y\,\nu({\rm d} y),
\end{align}
\begin{align}
c_2= \frac{1}{2}A +\frac{1}{2} \int_{\mathbb{R}\backslash \{0\}} y^2\,\nu({\rm d} y),
\end{align}
and
\begin{align}
c_k=  \int_{\mathbb{R}\backslash \{0\}}\frac{ (-y)^k}{k!} \,\nu({\rm d} y),\quad\quad k\ge 3.
\end{align}

\end{remark}

\section{Fokker-Planck equations for  SDEs with Marcus integrals}
\label{marcus888}

 In this section, we derive the Fokker-Planck equation for Marcus  SDE (\ref{Mainequation_strat}), in which  $L_t$ being a scalar L\'{e}vy process  with the triplet as  $(b, 1, \nu)$.
For simplicity,   we make the following assumption.

\textbf{Assumption ($H_2$)}\\
Assume that the noise intensity $\sigma(x,t)=\sigma(x) \neq 0$.

\medskip

Then (\ref{Mainequation_strat}) becomes
\begin{equation}\label{Main-equation}
dX_t=f(X_t, t)dt+\sigma(X_{t})\diamond dL_t.
\end{equation}

Define
\begin{align}\label{lamperti}
Y_t=H(X_t )=\int_0^{X_t}\frac{1}{\sigma(u )}du.
\end{align}
The transform $H$   in (\ref{lamperti}) is the Lamperti transform \cite{Iacus2008}.

By the chain rule of the Marcus integral \cite{Applebaum2009}, we have
\begin{eqnarray}\label{Additive-equation}
\nonumber {\rm d} Y_t&=&H'(X_{t})\diamond {\rm d}  X_t\\
\nonumber &=&H'(X_t)f(X_t, t){\rm d} t+\left(H'(X_t)\sigma(X_t)\right)\diamond {\rm d} L_t\\
\nonumber &=&\frac{f(X_t, t)}{\sigma(X_t)}{\rm d}t+{\rm d} L_t\\
&=&\frac{f(H^{-1}(Y_t), t)}{\sigma(H^{-1}(Y_t))}{\rm d} t+{\rm d} L_t.
\end{eqnarray}
Introducing
\begin{align}
\tilde f(Y_t, t)=\frac{f(H^{-1}(Y_t))}{\sigma(H^{-1}(Y_t))},
 \end{align}
equation (\ref{Additive-equation}) becomes
\begin{align}
 {\rm d} Y_t= \tilde f(Y_t, t) {\rm d} t +{\rm d} L_t.
\end{align}
By Theorem \ref{fpe-marcus} in the last section, the Forkker-Planck equation for $Y_t$ is
\begin{align}\label{additive_y}
\frac{\partial}{\partial t} p(y,t) &= -\frac{ \partial}{\partial y}\left (\tilde f(y,t)p(y,t)\right) -b \frac{\partial}{\partial y} p(y,t) +\frac{A}{2} \frac{\partial ^2}{\partial y^2} p(y,t)\nonumber\\
 &\quad + \int\limits_{{R \backslash \{0\}}} \left(p(y-r) -p(y)+ I_{(-1, 1)}(r)r\frac{\partial}{\partial y} p(y,t) \right) \nu(  {\rm d} r).
\end{align}

Let $q(x,t)$ represents the probability density function of $X_t$,   from the transform (\ref{lamperti}) and the Fokker-Planck equation (\ref{additive_y}), we get the following result for the desired Fokker-Planck equation for SDE (\ref{Main-equation}) defined by Marcus integrals.

\begin{theorem} \label{fpe-marcus}
Under the \textbf{Assumptions ($H_1$) and ($H_2$)}, the   Fokker-Planck equation for the Marcus SDE (\ref{Main-equation}) is
\small
\begin{align}\label{FPK_Marcus}
&\sigma(x) \frac{\partial}{\partial t} q(x,t) = -\sigma(x) \frac{ \partial}{\partial x}\left ( f(x,t)q(x,t)\right) -b \sigma(x) \frac{\partial}{\partial x}\left( \sigma(x)q(x,t)\right)  +\frac{A}{2}\sigma(x) \frac{\partial }{\partial x }\left(\sigma(x) \frac{\partial }{\partial x}\left(\sigma(x) q(x,t)\right)\right)\nonumber\\
&\quad
 + \int\limits_{{R \backslash \{0\}}} \left(\sigma(H^{-1} (H(x)-r))q(H^{-1} (H(x)-r), t) -\sigma(x)q(x,t) )+ I_{(-1,1)}(r)r \sigma(x)\frac{\partial}{\partial x} \left(\sigma(x)q(x,t)\right) \right) \nu (  {\rm d} r).
\end{align}
\normalsize
\end{theorem}

 \begin{remark}
 Consider a special case.
 When $\sigma(x)=1$, (\ref{FPK_Marcus}) reduces to
\begin{align}
\frac{\partial q(x,t)}{\partial t}
 &=-\frac{\partial }{\partial x}\left( f(x,t)q (x, t) \right) -b\frac{\partial }{\partial x}  q (x,t) ) + \frac{A}{2}  \frac{\partial^2}{\partial x^2} \ q(x,t) \nonumber\\
 &\quad + \int\limits_{R\backslash \{0\} }\left( q(x-y,t ) -q(x,t)+I_{(-1, 1)}(y)y \frac{\partial}{\partial x} q(x,t)  \right)\nu ({\rm d} y),
\end{align}
which is equivalent to (\ref{additive_FPK}), the Fokker-Planck equation for a SDE with an  additive L\'evy process. This is the consequence of the fact that for additive L\'{e}vy noise, the Marcus integral and It\^o integral are the same.

\end{remark}

\section{Examples}
\label{examples888}

In this section, we present the Fokker-Planck equations for some SDEs with special L\'evy processes, such as a Brownian motion together with  a Poisson process, a compound Poisson process, and finally a $\alpha$-stable L\'evy motion, respectively.

For simplicity, in all these examples, we take $f(x, t)=x$, and $\sigma(x,t)=x$. Then the It\^o SDE (\ref{sde_Ito}) and the Marcus SDE (\ref{Mainequation_strat}) become
\begin{align}\label{example1_sde1}
{\rm d}X_t= X_t {\rm d}t+X_t {\rm d} L_t,
\end{align}
and
\begin{align}\label{example_sde2}
{\rm d}X_t= X_t {\rm d}t+X_t \diamond{\rm d} L_t,
\end{align}
respectively.

{\bf{Example 1.}}
When  $L_t$  is a standard Brownian motion together with  a Poisson process with parameter $\lambda$,   the triplet  $(b, A, \nu)$ of $L_t$ is $(0, 1, \lambda \delta_1)$, where $\delta_1$ is the Dirac measure defined as
\begin{align}
\delta_1 (D) = \begin{cases}
1, &\text{if $1\in D$,}\\
0, &\text{if $1 \notin D$.}
\end{cases}
\end{align}
By using (\ref{multiplicative_FPK}), we get the Fokker-Planck equation for (\ref{example1_sde1})
\begin{align}
\frac{\partial p(x,t)}{\partial t}
 &=-\frac{\partial }{\partial x}\left(xp (x, t) \right) + \frac{1}{2}\frac{\partial ^2}{\partial x^2} \left(x^2 p(x,t)\right)+ \lambda \sum_{k=1}^\infty \frac{(-1)^k}{k!}  \frac{\partial ^k}{\partial x^k} \left(x^k p(x,t)\right).
\end{align}
By using (\ref{FPK_Marcus}), we get the Fokker-Planck equation for (\ref{example_sde2})
\begin{align}
x\frac{\partial q(x,t)}{\partial t} = -x  \frac{\partial  }{\partial x } \left(x q(x,t)\right) +  \frac{1}{2}x\frac{\partial }{\partial x } \left(x \frac{\partial }{\partial x }\left( xp(x,t)\right)\right) + \lambda e^{-1} |x|  q(e^{-1}|x|, t) - \lambda xq(x,t).
\end{align}

{\bf{Example 2.}}
When $L_t$ is a compound Poisson process defined as
\begin{align}
L_t=\sum_{i=1}^{N_t} \xi_i,
\end{align}
where $N_t$ is a Poisson process with parameter $\lambda$ and $\xi_i$ ($i=1, 2, \cdots$) are i.i.d. random numbers with standard normal distribution $\mathbf{N} (0,1)$,   the triplet $(b, A, \nu)$ of $L_t$  is $b=0$, $A=0$, and $\nu({\rm d}x)=\frac {\lambda }{\sqrt{2\pi}} e^{-\frac{x^2}{2}}{\rm d} x$. By using (\ref{multiplicative_FPK}), we get the Fokker-Planck equation for (\ref{example1_sde1})
\begin{align}
\frac{\partial p(x,t)}{\partial t}=-\frac{\partial }{\partial x}\left( f(x,t)p (x, t) \right) + \lambda \sum_{k=1}^{\infty}   c_{ k}  \frac{\partial ^{ k}}{\partial x^{ k}} \left(x^{ k}  p(x,t)\right),
\end{align}
where
\begin{align}
c_k=\begin{cases}
0, &\text{if $k$ is odd},\\
(k-1)!!, &\text{if $k$ is even}.
\end{cases}
&\text{\quad \quad ($k\in \mathbb{N}$),}
\end{align}
and $!!$ is the double factorial defined by
\begin{align}
n!!=\begin{cases}
n\cdot (n-2)  \cdots 1, &\text{if $n$ is odd},\\
n\cdot (n-2)\cdots 2, &\text{if $n$ is even}.
\end{cases}
\end{align}
By using (\ref{FPK_Marcus}), we get the Fokker-Planck equation for (\ref{example_sde2})
\begin{align}
x\frac{\partial q(x,t)}{\partial t} = -x  \frac{\partial  }{\partial x } \left(x q(x,t)\right) + \lambda \int_{\mathbb{R}\backslash \{0\}} \left(  e ^{-r}|x| q(e ^{-r}|x|, t) - xq(x,t) \right) {\rm d} r.
\end{align}

{\bf{Example 3.}}
When $L_t$ is a symmetric $\alpha$-stable L\'evy motion with the triplet as $b=1$, $A=0$ and $\nu({\rm d}x)= \frac{ {C_\alpha \rm d} x}{|x|^{1+\alpha}}$, where $C_\alpha$ is a constant (depending on $\alpha$). Then by using (\ref{multiplicative_FPK}), we get the Fokker-Planck equation for (\ref{example1_sde1})
\begin{align}
&\frac{\partial p(x,t)}{\partial t}
 =-\frac{\partial }{\partial x}\left( f(x,t)p (x, t) \right)\nonumber\\
 & \quad +C_\alpha \int\limits_{\mathbb{R}\backslash \{0\} }  \left[ \sum_{k=1}^{\infty} \frac{(-y)^k}{k!} \frac{\partial ^k}{\partial x^k} \left(\sigma^k (x,t) p(x,t)\right) + I_{(-1, 1)}(y) y  \frac{\partial}{\partial x} \left(\sigma (x,t)p(x,t)\right) \right]  \, \frac{ { \rm d} y}{|y|^{1+\alpha}}.
\end{align}
By using (\ref{FPK_Marcus}), we get the Fokker-Planck equation for (\ref{example_sde2})
\begin{align}
&x\frac{\partial q(x,t)}{\partial t} = -x  \frac{\partial}{\partial x } \left(x q(x,t)\right) \nonumber\\
& \quad + \lambda C _\alpha\int_{\mathbb{R}\backslash \{0\}} \left[  e ^{-r}|x| q(e ^{-r}|x|, t) - xq(x,t)    I_{(-1, 1)}(r) r  \frac{\partial}{\partial x} \left(\sigma (x,t)p(x,t)\right) \right]  \, \frac{ { \rm d} r}{|r|^{1+\alpha}}.
\end{align}

\newpage

\end{document}